\documentclass[11pt]{article}

\usepackage{amsmath,amsfonts,amsthm,amssymb,graphicx,tikz}
\usepackage[all,2cell,ps]{xy}

\theoremstyle{plain}
\newtheorem{thm}{Theorem}[section]

\newtheorem{lem}[thm]{Lemma}

\newtheorem{cor}[thm]{Corollary}

\theoremstyle{definition}

\newtheorem{rem}[thm]{Remark}

\theoremstyle{remark}

\newcommand{\bbA}{\mathbb{A}}
\newcommand{\bbB}{\mathbb{B}}
\newcommand{\bbC}{\mathbb{C}}

\newcommand{\bbQ}{\mathbb{Q}}
\newcommand{\bbR}{\mathbb{R}}

\newcommand{\bbZ}{\mathbb{Z}}

\newcommand{\calE}{\mathcal{E}}

\newcommand{\calH}{\mathcal{H}}

\newcommand{\calJ}{\mathcal{J}}

\newcommand{\calO}{\mathcal{O}}
\newcommand{\calP}{\mathcal{P}}
\newcommand{\calQ}{\mathcal{Q}}

\newcommand{\fraks}{\mathfrak{s}}

\newcommand{\al}{\alpha}

\newcommand{\Gam}{\Gamma}

\newcommand{\de}{\delta}
\newcommand{\Del}{\Delta}

\newcommand{\Lam}{\Lambda}
\newcommand{\sig}{\sigma}

\DeclareMathOperator{\M}{M}
\DeclareMathOperator{\SL}{SL}
\DeclareMathOperator{\PSL}{PSL}

\DeclareMathOperator{\PU}{PU}
\DeclareMathOperator{\SU}{SU}

\DeclareMathOperator{\Aut}{Aut}

\DeclareMathOperator{\Gal}{Gal}

\DeclareMathOperator{\Res}{Res}

\newcommand{\conj}{\overline}

\newenvironment{pf}{\begin{proof}}{\end{proof}}

\newenvironment{enum}{\begin{enumerate}}{\end{enumerate}}

\makeatletter
\let\@@pmod\pmod
\DeclareRobustCommand{\pmod}{\@ifstar\@pmods\@@pmod}
\def\@pmods#1{\mkern4mu({\operator@font mod}\mkern 6mu#1)}
\makeatother

\usetikzlibrary{arrows}

\title{Lattices in $\PU(n,1)$ that are not profinitely rigid}
\author{Matthew Stover\footnote{This material is based upon work supported by Grant Number 523197 from the Simons Foundation/SFARI.} \\ \small{Temple University}\\ \small{\textsf{mstover@temple.edu}}}
\date{\today}

\begin{document}

\maketitle

\begin{abstract}
Using conjugation of Shimura varieties, we produce nonisomorphic, cocompact, torsion-free lattices in $\PU(n,1)$ with isomorphic profinite completions for all $n \ge 2$. This disproves a conjecture of D.\ Kazhdan and gives the first examples nonisomorphic lattices in a semisimple Lie group of real rank one with isomorphic profinite completions, answering two questions of A.\ Reid.
\end{abstract}

\section{Introduction}\label{sec:Intro}

Let $\Gam$ be a finitely generated residually finite group. Its profinite completion $\hat{\Gam}$ is a powerful tool in studying not only the group-theoretic properties of $\Gam$, but also the geometric properties of a manifold or variety with topological fundamental group $\Gam$. Recently, this connection has led to a great deal of interest in \emph{profinite rigidity}; one says that $\Gam$ is profinitely rigid if any finitely generated residually finite group $\Del$ with $\hat{\Gam} \cong \hat{\Del}$ is necessarily isomorphic to $\Gam$. See Reid's ICM address \cite{ReidICM} for an introduction to this topic and \cite{ReidStA} for a more expanded survey.

Profinite rigidity is of particular interest for the groups one encounters in geometric and low-dimensional topology. For example, Bridson, Conder, and Reid proved that if two finitely generated Fuchsian groups (i.e., discrete subgroups of $\PSL_2(\bbR)$) have isomorphic profinite completions, then the two Fuchsian groups are necessarily isomorphic \cite{BCR}. More recently, Bridson, McReynolds, Reid, and Spitler \cite{BMRS} showed that the arithmetic Kleinian group $\PSL_2(\bbZ[e^{2 \pi i / 3}])$ is profinitely rigid. In the negative direction, Baumslag gave examples of nonisomorphic infinite finitely generated nilpotent groups with isomorphic profinite completions \cite{Baumslag}, and the congruence subgroup property leads to examples of nonisomorphic higher rank lattices in semisimple Lie groups with isomorphic profinite completions \cite{Aka}.

Motivated by the above, Reid asked whether or not lattices in rank one real Lie groups are profinitely rigid \cite[Q.\ 10]{ReidStA}. The purpose of this paper is to answer this question in the negative, confirming Reid's speculation.

\begin{thm}\label{thm:Main}
For each $n \ge 2$, there exists a pair of nonisomorphic, cocompact, torsion-free lattices $\Gam, \Del < \PU(n,1)$ such that $\hat{\Gam} \cong \hat{\Del}$. One can choose $\Gam$ and $\Del$ to be noncommensurable.
\end{thm}

We also note that one can guarantee that these lattices do not have the congruence subgroup property. In fact, we will show that they can have arbitrarily large first betti number; see Lemma \ref{lem:b1}. On the other hand, our methods will always produce examples with the same volume; see Lemma \ref{lem:Vol} and compare with \cite[Q.\ 7.4]{ReidICM}.  We can also produce arbitrarily long sequences of nonisomorphic lattices with isomorphic profinite completions. Since cocompact lattices in $\PU(n,1)$ are Gromov hyperbolic groups, we also obtain the following, answering another question of Reid \cite[Q.\ 14]{ReidStA}.

\begin{thm}\label{thm:Gromov}
For each $d \ge 2$, there exist pairwise-nonisomorphic torsion-free Gromov hyperbolic groups $\Gam_1, \dots, \Gam_d$ such that $\hat{\Gam}_j \cong \hat{\Gam}_k$ for all pairs $1 \le j,k \le d$.
\end{thm}

Our proof uses a combination of techniques from discrete subgroups of Lie groups and algebraic geometry. That one can find distinct smooth projective surfaces with nonisomorphic fundamental groups having the same profinite completion is well-known, and goes back to Serre \cite{Serre}, who noticed that it suffices to construct complex varieties $V, W$ with nonisomorphic fundamental groups such that $W = \tau V$ for some $\tau \in \Aut(\bbC)$. It follows immediately that they have the same \'etale fundamental group, i.e., that their topological fundamental groups have isomorphic profinite completions.

We employ this basic idea via conjugation of Shimura varieties. Results of this kind date back to work of Kazhdan in the early 70s \cite{KazhdanICM}, and we closely follow more recent work of Milne and Suh \cite{MilneSuh}, who constructed examples of nonisomorphic higher-rank lattices with isomorphic profinite completions using conjugate but not homeomorphic Shimura varieties. (In contrast with \cite{Aka}, they do not need to assume the congruence subgroup property holds.) The main contribution of this paper is to show that their use of the Margulis superrigidity theorem can be weakened to Mostow rigidity, which then applies for $\PU(n,1)$.

Said differently, we use the theory of arithmetic groups to find smooth projective surfaces $V, W$ with universal cover $\bbB^n$ such that $\pi_1(V) \not\cong \pi_1(W)$ but $W = \tau V$ for some $\tau \in \Aut(\bbC)$. Here $\bbB^n$ denotes the unit ball in $\bbC^n$ with its (negatively curved) hermitian symmetric Bergman metric. In particular, we will show the following, which disproves a conjecture of Kazhdan \cite[p.\ 325]{KazhdanICM}.

\begin{thm}\label{thm:KazhdanNo}
For each $n \ge 2$, there exists a pair of smooth projective varieties $V, W$ defined over a number field $F$ such that
\begin{enum}

\item $V$ and $W$ have nonisomorphic (topological) fundamental groups;

\item $V$ and $W$ both have universal covering biholomorphic to the ball $\bbB^n$ with its Bergman metric;

\item there exists $\tau \in \Aut(\bbC / F)$ such that $W = \tau V$.

\end{enum}
\end{thm}

We now briefly discuss the organization of the paper. In \S \ref{sec:Prelim}, we give preliminary material from Milne--Suh \cite{MilneSuh} and prove a generalization of their main result using Mostow rigidity. In \S \ref{sec:Ex}, we give very concrete examples of lattices in $\PU(2,1)$ for which the results in \S \ref{sec:Prelim} will apply. These examples suffice to prove the theorems discussed in this introduction.

\subsubsection*{Acknowledgments} The author thanks Alan Reid for a number of stimulating conversations related to the topic of this paper and Henry Wilton for asking about volume.

\section{Preliminaries}\label{sec:Prelim}

The purpose of this section is to translate the results of Milne and Suh \cite{MilneSuh} to the situation of interest for us, namely the unit ball $\bbB^n$ in $\bbC^n$ with its Bergman metric. We assume familiarity with the basic notions used there. This metric makes $\bbB^n$ into a Hermitian symmetric domain with $\PU(n,1)$ the identity component of the group $\mathrm{Hol}^+(\bbB^n)$ of holomorphic automorphisms of $\bbB^n$.

To simplify the discussion, we will only consider the arithmetic subgroups of $\PU(n,1)$ of `simplest type', though our methods apply to an arbitrary arithmetic lattice. Let $F$ be a totally real number field with degree $d$ over $\bbQ$ and $\nu_1, \dots, \nu_d$ be the distinct embeddings of $F$ into $\bbR$. Fix a totally imaginary quadratic extension $E / F$ and let $\sig$ denote the nontrivial Galois involution for the extension. Choose a $\sig$-hermitian form $h$ on $E^{n+1}$ and let $H$ be the $F$-algebraic group $\SU(h)$ of special unitary automorphisms of $h$. This is an absolutely almost simple, simply connected $F$-algebraic group.

For each embedding $\nu_j : F \to \bbR$ we obtain a real algebraic group $H_{\nu_j}$ with $H_{\nu_j}(\bbR) = H(F \otimes_{\nu_j} \bbR)$, which is then isomorphic to $\SU(p_j,q_j)$, where $p_j+q_j = n+1$ and $(p_j,q_j)$ is the signature of $h$ considered as a hermitian form on $\bbC^{n+1}$ via this embedding. Our standing assumption in what follows is that there is exactly one $1 \le j \le d$ such that $H_{\nu_j}(\bbR) \cong \SU(n,1)$ and that $H_{\nu_k}(\bbR) \cong \SU(n+1)$ for all $k \neq j$.\footnote{It is tradition to rearrange the places of $F$ so that $\nu_1$ is the place where the signature is $(n,1)$. For the arguments in this paper, this assumption will only complicate things.} The restriction of scalars $H_* = \Res_{F / \bbQ}(H)$ is then a $\bbQ$-algebraic group with
\[
H_*(\bbR) \cong \SU(n, 1) \times \SU(n+1)^{d-1},
\]
and hence $H_*(\bbR)$ admits a projection onto $\PU(n,1)$ with compact kernel.

Let $\bbA_F^\infty$ be the finite adeles of $F$ and $K$ an open compact subgroup of $H(\bbA_F^\infty)$. The group $\Lam = H(F) \cap K$ is a \emph{congruence subgroup} of $H(F)$, and its image $\Gam$ in $\PU(n,1)$ is a lattice, i.e., a discrete subgroup so that $V = \bbB^n / \Gam$ has finite volume. This class of lattices includes the more traditional congruence subgroups defined by matrices congruent to the identity modulo some ideal; see \S \ref{sec:Ex} for an explicit example connecting the two notions. It is then a basic fact that $V$ is compact if and only if $d > 1$, i.e., $F \neq \bbQ$. Moreover, $V$ is a normal quasiprojective variety, hence it is projective when $V$ is compact. For sufficiently small $K$ we can assume that $V$ is smooth, which is equivalent to $\Gam$ being torsion-free; in the compact case $V$ is then a smooth complex projective variety of dimension $n$.

In the language of \cite{MilneSuh}, such a $V$ is a Shimura variety of type $(H, \bbB^2)$. The theory of Shimura varieties then implies that $V$ has a model defined over a number field called its \emph{reflex field}. Given any element $\tau$ of $\Aut(\bbC)$, we then apply $\tau$ to the polynomials defining $V$ as a smooth quasiprojective variety and obtain a possibly new variety $\tau V$. It is a fundamental result going back to Kazhdan \cite{KazhdanICM} that $\tau V$ is also a Shimura variety.\footnote{In fact, Kazhdan proved the the same holds when $V$ is an arbitrary arithmetic quotient of a hermitian symmetric domain, not just the quotient by a congruence subgroup.} Milne and Suh gave the following precise formulation of this:

\begin{thm}[Thm.\ 1.3 \cite{MilneSuh}]\label{thm:MilneSuh}
Let $V$ be a smooth algebraic variety over $\bbC$, and let $\tau$ be an automorphism of $\bbC$. If $V$ is of type $(H, X)$ with $H$ simply connected, then $\tau V$ is of type $(H^\prime, X^\prime)$ for some semisimple algebraic group $H^\prime$ over $F$ such that:
\[
\begin{cases}
H^\prime_\nu \cong H_{\tau \circ \nu}\ \textrm{for all real embeddings}\ \nu : F \to \bbR,\, \textrm{and} \\
H^\prime_\nu \cong H_\nu\ \textrm{for all nonarchimedean places}\ \nu\ \textrm{of}\ F
\end{cases}
\]
\end{thm}

We now prove a result that expands the scope of the examples that one can produce using the techniques of \cite{MilneSuh}. The key observation is that when $X$ is not the Poincar\'e disk and $\Gam$ is an irreducible lattice in $\mathrm{Hol}^+(X)$, one can replace the use of Margulis superrigidity in \cite{MilneSuh} with Mostow rigidity and obtain the same conclusions.

\begin{thm}\label{thm:ApplyMostow}
Suppose that $V$, $H$, $\tau$, and $H^\prime$ are as in Theorem \ref{thm:MilneSuh}. Furthermore, assume that $\dim_\bbC(V) > 1$ and that $H$ is almost simple over $F$. If $H^\prime$ is not isomorphic as an $F$-algebraic group to $\sig H$ for any $\sig \in \Aut(F / \bbQ)$, then $\pi_1(V)$ is not isomorphic to $\pi_1(\tau V)$.
\end{thm}

\begin{pf}
We have that $V$ is a Shimura variety of type $(H, X)$ and $\tau V$ has type $(H^\prime, X^\prime)$. The assumption that $H^\prime$ is not isomorphic to $\sig H$ for any $\sig$ in $\Aut(F / \bbQ)$ is equivalent to the assumption that $H_* = \Res_{F / \bbQ}(H)$ and $H^\prime_* = \Res_{F / \bbQ}(H^\prime)$ are not isomorphic as $\bbQ$-algebraic groups \cite[\S 3]{BorelDense}. Indeed, such an isomorphism of $\bbQ$ algebraic groups is necessarily the composition of a $\bbQ$-algebra isomorphism $\sig : F \to F$ with an isomorphism $\sig H \to H^\prime$. That $H$ is almost simple over $F$ implies that $\Gam = \pi_1(V)$ is an irreducible lattice in $\mathrm{Hol}^+(X)$. Let $\Gam^\prime = \pi_1(\tau V)$.

Suppose that $\Gam \cong \Gam^\prime$. Since $\dim_\bbC(V) > 1$, Mostow rigidity holds. Therefore $X = X^\prime$ and $V$ is isometric to $\tau V$, hence there is some $g \in \mathrm{Hol}^+(X)$ such that $\Gam^\prime = g \Gam g^{-1}$. This conjugation defines an isomorphism between the commensurators $\mathrm{C}_{\textrm{Hol}}(\Gam)$ and $\mathrm{C}_{\textrm{Hol}}(\Gam^\prime)$ of $\Gam$ and $\Gam^\prime$ in $\mathrm{Hol}^+(X)$, respectively. It is a fundamental result due to Borel \cite[Thm.\ 2]{BorelDense} that $\mathrm{C}_{\textrm{Hol}}(\Gam)$ is isomorphic to $\conj{H}_*(\bbQ)$, where $\conj{H}_*$ is the adjoint form of the $\bbQ$-algebraic group $H_*$, and similarly for $\mathrm{C}_{\textrm{Hol}}(\Gam^\prime)$ and $\conj{H}^\prime_*$. It follows that $\conj{H}_*$ and $\conj{H}^\prime_*$ are isomorphic as $\bbQ$-algebraic groups, and hence $H_*$ is similarly isomorphic to $H^\prime_*$. This contradicts our assumption that $H^\prime$ is not isomorphic to $\sig H$ for any $\sig \in \Aut(F / \bbQ)$ and completes the proof of the theorem.
\end{pf}

\begin{rem}
As noticed in \cite[Rem.\ 2.7]{MilneSuh}, it might still be the case that $\pi_1(V)$ is not isomorphic to $\pi_1(\tau V)$ when $H \cong H^\prime$. Indeed, the two fundamental groups might only be commensurable but not isomorphic. It would be interesting to find such an example.
\end{rem}

\section{An family of examples}\label{sec:Ex}

To justify the main results in the introduction, it suffices to produce examples of algebraic groups $H$ to which Theorem \ref{thm:ApplyMostow} applies. We do this for $\PU(2,1)$ for concreteness, though our examples trivially generalize to $\PU(n,1)$ for any $n \ge 3$. We describe in detail the construction of two nonisomorphic torsion-free cocompact lattices in $\PU(2,1)$ with the same profinite completion. At the end of the section, we indicate how our examples immediately generalize to give arbitrarily large families. Very similar higher-rank examples appear in \cite{MilneSuh}.

\medskip

Consider the totally real cubic $F = \bbQ(\al)$ of discriminant $148$ where $\al$ has minimal polynomial
\[
p(x) = x^3 - x^2 - 3 x + 1.
\]
Then $F$ has Galois group the symmetric group $S_3$ on three letters and $\Aut(F / \bbQ)$ is trivial. There is exactly one real embedding $\nu_1$ of $F$ with $\nu_1(\al) > 0$, and there are two embeddings $\nu_2, \nu_3$ with $\nu_j(\al) < 0$.

The Galois closure $L$ of $F$ is totally real with minimal polynomial $q(x) = x^6 - 20 x^4 + 100 x^2 - 148$ and discriminant $810448$. Write $L = \bbQ(\beta)$, where $\beta$ is a root of $q(x)$. The three distinct embeddings of $F$ into $L$ are given by
\begin{align*}
\rho_1(\al) &= \frac{3}{2} \beta^4 - 25 \beta^2 + 67 \\
\rho_2(\al) &= -\frac{3}{4} \beta^4 + \frac{25}{2} \beta^2 + \frac{1}{2} \beta - 33 \\
\rho_3(\al) &= -\frac{3}{4} \beta^4 + \frac{25}{2} \beta^2 - \frac{1}{2} \beta - 33
\end{align*}
We identify $S_3 = \Gal(L / \bbQ)$ with group of permutations of the set $\{\nu_1, \nu_2, \nu_3\}$ of real embeddings of $F$.

Let $E / F$ be any totally imaginary quadratic extension. For concreteness, we take $E = F(\de)$ with $\de^2 = -1$. Let $\rho$ denote the nontrivial Galois automorphism of $E / F$. We fix the two $\rho$-hermitian forms
\begin{align*}
h_1 &= \begin{pmatrix} -\al & 0 & 0 \\ 0 & -\al & 0 \\ 0 & 0 & -1 \end{pmatrix} \\
h_2 &= \begin{pmatrix} 2-\al^2 & 0 & 0 \\ 0 & 2-\al^2 & 0 \\ 0 & 0 & -1 \end{pmatrix}
\end{align*}
on $E^3$, and let $H_j = \SU(h_j)$ be the special unitary group of $h_j$. Each $H_j$ is an absolutely almost simple, simply connected $F$-algebraic group.

Then we have:
\begin{align*}
\Res_{F / \bbQ}(H_1)(\bbR) &\cong \SU(2,1) \times \SU(3) \times \SU(3) \\
\Res_{F / \bbQ}(H_2)(\bbR) &\cong \SU(3) \times \SU(2,1) \times \SU(3)
\end{align*}
In the first case, $\nu_j(h_1) \in \M_3(F \otimes_{\nu_j} \bbR)$ determines a signature $(2,1)$ hermitian form on $\bbC^3$ for $j = 1$ and definite form for $j = 2,3$. For the second, the same is true with the roles of $\nu_1$ and $\nu_2$ exchanged.

\begin{lem}\label{lem:SamePadic}
Let $\nu$ be a nonarchimedean place of $F$ and $F_\nu$ the associated local field. Then $H_1(F_\nu) \cong H_2(F_\nu)$.
\end{lem}

\begin{pf}
When $\nu$ splits in $E$, one has $H_j(F_\nu) \cong \SL_3(F_\nu)$. When $\nu$ is inert or ramifies and $\mu$ is the unique prime of $E$ lying over $\nu$, $H_j(F_\nu)$ is isomorphic to the unique special unitary group in $3$ variables over $F_\nu$ with respect to the quadratic extension $E_\mu / F_\nu$. See \cite[Ch.\ 10]{Scharlau} for details on hermitian forms over global fields and their localizations.
\end{pf}

This immediately implies the following.

\begin{cor}\label{cor:AdeleIso}
Let $\bbA_F^\infty$ be the finite adeles of $F$. Then there is an isomorphism
\[
H_1(\bbA_F^\infty) \cong H_2(\bbA_F^\infty).
\]
\end{cor}

Fix an open compact subgroup $K < H_1(\bbA_F^\infty) \cong H_2(\bbA_F^\infty)$ and define
\[
\Lam_j = H_j(F) \cap K,
\]
where $\Lam_2$ is determined by a fixed isomorphism $H_1(\bbA_F^\infty) \to H_2(\bbA_F^\infty)$. Then each $\Lam_j$ defines a lattice $\Gam_j < \PU(2,1)$ via the real embedding $\nu_j$ of $F$. Since $H_j$ is anisotropic, it follows that $\Gam_j$ is cocompact. We assume that $K$ is sufficiently small that each $\Gam_j$ is torsion-free.

For a concrete example, for any prime ideal $\calP$ of the ring of integers $\calO_F$ of $F$, we can choose $K$ such that $\Lam_1$ is the congruence subgroup of $H_1(\calO_F)$ consisting of those matrices in $\SL_3(\calO_E)$ that preserve $h_1$ and are congruent to the identity modulo $\calP$. Let $\nu_\calQ$ denote the place of $F$ associated with the prime ideal $\calQ$ of $\calO_F$ and $\calO_\calQ \subset F_{\nu_\calQ}$ be the integral closure of $\calO_F$ in $F_{\nu_\calQ}$. We then have that this classical congruence subgroup is associated with the compact open subgroup
\[
K_\calP = H_1(\calP \calO_\calP) \times \prod_{\calQ \neq \calP} H_1(\calO_\calQ) < H_1(\bbA_F^\infty).
\]
For all but finitely many $\calP$, this determines a torsion-free lattice, and this generalizes in the obvious way to any ideal $\calJ$ of $\calO_F$.

However, notice that this does not necessarily mean that $\Lam_2$ has the same description as a congruence subgroup in the more well-known sense, since the identification of $K$ with an open compact subgroup of $H_2(\bbA_F^\infty)$ depends on the isomorphisms between $H_1(F_\nu)$ and $H_2(F_\nu)$ for each nonarchimedean place $\nu$ of $F$. (We will not use this, but remark that $\al$, $\al^2-2$, and $-1$ generate the unit group of $F$; this could be used to explicitly calculate the image of $K$ in $H_2(\bbA_F^\infty)$.)

\begin{thm}\label{thm:SamePro}
The lattices $\Gam_1, \Gam_2 < \PU(2,1)$ described above are not isomorphic but have isomorphic profinite completions.
\end{thm}

\begin{pf}
Let $V_j$ be the smooth complex hyperbolic manifold $\bbB^2 / \Gam_j$, considered as a smooth projective surface defined over a number field. We will show that there is an element $\tau \in \Aut(\bbC)$ such that $V_2 = \tau V_1$. It follows that the two varieties have the same \'etale fundamental groups. For example, see the discussion in \cite[\S 1.1]{MilneSuh}. Recall that the \'etale fundamental group of $V_j$ is the profinite completion of the topological fundamental group, which is $\Gam_j$. Since there is by construction no element $\sig$ in $\Aut(F / \bbQ) = \{1\}$ such that $\sig H_1$ is $F$-isomorphic to $H_2$, it follows from Theorem \ref{thm:ApplyMostow} that $\Gam_1$ and $\Gam_2$ are not commensurable in $\PU(2,1)$, hence they certainly are not isomorphic.

Let $\tau$ be an element of $\Aut(\bbC)$ whose restriction to $L$ is the transposition $(1\, 2) \in S_3$, where $S_3$ is identified with the permutation group on the set $\{\nu_1, \nu_2, \nu_3\}$ of real embeddings of $F$. Then Theorem \ref{thm:MilneSuh} implies that $\tau V_1$ is isomorphic to a Shimura variety $V_1^\tau$ associated with an $F$-algebraic group $H_1^\prime$ with the property that $H_1^\prime(F_\nu) \cong H_1(F_\nu)$ for every nonarchimedean place $\nu$ of $F$ and
\[
H_1^\prime(F_\nu) \cong H_1(F_{\tau(\nu)})
\]
for every real embedding. By construction we then have that $H_1^\prime(F_\nu) \cong H_2(F_\nu)$ for every place $\nu$ of $F$, so the Hasse principle for unitary groups \cite[p.\ 359]{PandR} implies that $H_1^\prime \cong H_2$ as $F$-algebraic groups. Fixing the isomorphism between $H_1(\bbA_F^\infty)$ and $H_2(\bbA_F^\infty)$ as in \cite[Rem.\ 1.6]{MilneSuh} to define $\Lam_2$ implies that $\pi_1(\tau V_1) \cong \Gam_2$. This completes the proof of the theorem.
\end{pf}

\begin{rem}
To prove the theorem, we did not need to prove that $H_1^\prime$ is isomorphic to $H_2$, only that $H_1^\prime$ is not $F$-isomorphic to $\sig H_1$ for some $\sig \in \Aut(F / \bbQ)$, which is obvious from the choice of $\tau$. We included $H_2$ to make the identification more concrete.
\end{rem}

We now indicate how one bootstraps the above to get arbitrarily long sequences $\Gam_1, \dots, \Gam_d$ of torsion-free cocompact lattices in $\PU(2,1)$ with isomorphic profinite completions. As in \cite{MilneSuh}, one replaces $F$ with a totally real field of degree $d$ with trivial automorphism group. One then finds $\Aut(F / \bbQ)$-inequivalent algebraic groups $H_1, \dots, H_d$ for which $H_j = \tau_j H_1$ for some $\tau_j \in \Aut(\bbC)$ and implements the above strategy to construct $V_j = V_1^{\tau_j}$. (Note that the above actually allows one to find three nonisomorphic lattices by choosing $h_3$ that has signature $(2,1)$ at $\nu_3$.)

\medskip

\noindent
We close the paper by proving some of the features of our examples indicated in the introduction.

\begin{lem}\label{lem:b1}
One can choose the lattices $\Gam_1, \Gam_2$ in Theorem \ref{thm:Main} to have arbitrarily large first betti number. In particular, they do not have the congruence subgroup property.
\end{lem}

\begin{pf}
It is well-known that having a homomorphism onto $\bbZ$ implies that a lattice does not have the congruence subgroup property. See \cite[Thm.\ 7.2]{LubotzkySegal}. For the lattices of `simplest type' considered in this paper, Kazhdan \cite{KazhdanB1} and Shimura \cite{ShimuraB1} showed that we can choose the congruence subgroup of $H_1(F)$ to have positive first betti number. It is then an observation due to Borel \cite{BorelB1} that one can find a congruence subgroup for which the first betti number is arbitrarily large.
\end{pf}

\begin{rem}
Since the groups $\Gam_1$ and $\Gam_2$ have the same profinite completions, they have the same set of finite abelian quotients. In particular, their abelianizations are isomorphic.
\end{rem}

\begin{lem}\label{lem:Vol}
The lattices $\Gam_1$ and $\Gam_2$ have the same covolume, i.e., $\bbB^n / \Gam_1$ and $\bbB^n / \Gam_2$ have the same volume.
\end{lem}

\begin{pf}
This is a consequence of Prasad's volume formula \cite[Thm.\ 3.7]{Prasad}. With notation as above, we have
\[
\Lam_j = H_j(F) \cap K,
\]
for $K < H_j(\bbA_F^\infty)$ an open compact subgroup, so $\Lam_j$ is a lattice in $\SU(n, 1)$ and $\Gam_j$ is its projection to $\PU(n,1)$. In particular, recall that $K$ is independent of $j = 1,2$. All the terms in Prasad's formula for the covolume of $\Lam_j$ depend on:
\begin{enum}

\item the totally real field $F$,

\item the imaginary quadratic extension $E$,

\item the dimension of $H_j$,

\item a number $\fraks(\calH_j)$ depending only on the unique absolutely quasi-simple, simply connected group $\calH_j$ defined and quasi-split over $F$ \cite[\S 0.4]{Prasad},

\item a product depending on the `exponents' $\{m_i\}$ of $\calH_j$,

\item the Tamagawa number of $H_j$, which is $1$ \cite[\S 3.3]{Prasad},

\item and an Euler factor $\calE_j$.

\end{enum}
Since $\dim(H_1) = \dim(H_2)$ and $\calH_1 \cong \calH_2$ is the unique unitary group with respect to $E / F$ that is everywhere quasi-split, the factors in (1) - (6) above are all the same. Critically, the Euler factors $\calE_j$ for $\Lam_j$ are completely determined by the open compact subgroup $K$, which is independent of $j$ (more specifically, it is determined by the parahoric subgroups $K_\nu = H_j(F_\nu) \cap K$ for each nonarchimedean place $\nu$ of $F$). This proves that $\Lam_1$ and $\Lam_2$ have the same covolume in the sense of \cite{Prasad}. The volume of $\bbB^n / \Gam_j$ is a constant multiple of the covolume of $\Lam_j$ (e.g., see \cite[\S 4.2]{BorelPrasad}), so the lemma follows.
\end{pf}

\bibliography{ProfiniteFlabby}

\end{document}